\documentclass[10pt,twoside]{article}
\usepackage{amssymb}
\usepackage[mathscr]{eucal}
\usepackage{eufrak}
\usepackage{amsmath,amsthm}
\usepackage{mathrsfs}
\usepackage[colorlinks=true,
   urlcolor=blue,           filecolor=green,      
   citecolor=green,      
   linkcolor=red,           bookmarks=true,
  unicode,
   plainpages=false,   ]{hyperref}
\usepackage{color}
\usepackage[latin1]{inputenc}
\def\bp{\begin{proof}}
\def\ep{\end{proof}}

\def\n{\nabla}
\def\intl#1{\int\limits_{#1}}
\def\intll#1#2{\int\limits_{#1}^{#2}}
\def\dm{|\hskip-0.05cm|}
\def\OO{\Omega}
\def\displ{\displaystyle}
\def\VSE{\vspace{6pt}\\&\displ }
\def\VS{\vspace{6pt}\\\displ }
\def\rf#1{{\rm(\ref{#1})}}
\def\chiu{\hfill$\displaystyle\vspace{4pt}
\underset\Box\null$\par}
\def\R{\Bbb R}
\def\N{\Bbb N}

\def\à{\`{a}}

\def\vep{\varepsilon}

\def\be{\begin{equation}}
\def\ba{\begin{array}}
\def\ea{\end{array}}
\def\ee{\end{equation}}
\def\vs1{\vspace{1ex}}

\def\ov{\overline}
\def\po{\partial\Omega}

\def\é{\'{e}}
\font\sc=cmcsc10
\title{\large\bf Navier-Stokes equations:\\ an analysis of a possible gap to achieve the energy equality}
\author{\sc  Francesca Crispo\thanks{Dipartimento di Matematica e Fisica,  
Universit\`{a} degli Studi della Campania
``L. Vanvitelli'', via Vivaldi 43, 81100 Caserta,
 Italy.
 francesca.crispo@unicampania.it}
\and \sc
Carlo Romano Grisanti \thanks{Dipartimento di Matematica, Universit\`a di Pisa, via Buonarroti 1/c, 56127 Pisa, Italy. carlo.romano.grisanti@unipi.it}\and \sc  Paolo Maremonti
\thanks{
Dipartimento di Matematica e Fisica,  
Universit\`{a} degli 
Studi della Campania
``L. Vanvitelli'', via Vivaldi 43, 81100 Caserta,
 Italy.
paolo.maremonti@unicampania.it}\thanks{ The research activity of F. Crispo and P. Maremonti is performed under the
auspices of   GNFM-INdAM. The research activity of F. Crispo is  also supported by the Program  ``Vanvitelli per la Ricerca: VALERE 2019'' financed by the University of Campania ``L. Vanvitelli''. The research activity of C.R. Grisanti is performed under the auspices
of  GNAMPA-INdAM and partially supported by the Research
Project of the University of Pisa "Energia e regolarità: tecniche
innovative per problemi classici di equazioni alle derivate parziali".}}
\date{}
\begin{document}
\markboth{\footnotesize\rm Francesca Crispo, Carlo Romano Grisanti and   P.
Maremonti} {\footnotesize\rm
Navier-Stokes equations: an analysis of a possible gap to achieve the energy equality}
\maketitle 
{\bf Abstract} - {\small The paper is concerned with 
the IBVP of  the Navier-Stokes equations. The goal is the attempt to construct a weak solution enjoying  an energy equality. This result is a natural continuation and improvement of the one obtained by the same authors in \cite{CGM}. } 
 \par\noindent
 \vskip -0.7true cm\noindent
\newcommand{\red}{\protect\bf}
\renewcommand\refname{\centerline
{\red {\normalsize \bf References}}}
\newtheorem{ass}
{\bf Assumption} 
\newtheorem{defi}
{\bf Definition} 
\newtheorem{tho}
{\bf Theorem} 
\newtheorem{rem}
{\sc Remark} 
\newtheorem{lemma}
{\bf Lemma} 
\newtheorem{coro}
{\bf Corollary} 
\newtheorem{prop}
{\bf Proposition} 
\renewcommand{\theequation}{\arabic{equation}}
\setcounter{section}{0}
\section{Introduction}\label{intro} This note concerns the 3D-Navier-Stokes initial boundary 
value problem:
\be\label{NS}\ba{l}v_t+v\cdot
\nabla v+\nabla\pi_v=\Delta
v,\;\nabla\cdot
v=0,\mbox{ in }(0,T)\times\OO,\VS v=0\mbox{ on }(0,T)\times\po,\hskip0.12cm
v(0,x)=v_0(x)\mbox{ on
}\{0\}\times\OO.\ea\ee  In system \rf{NS} $\OO\subseteq\R^3$ is assumed to be  bounded or exterior, and its boundary is  smooth. The symbol   $v$ denotes the kinetic field, $\pi_v$ is the pressure field,
 $v_t:=
\frac\partial{\partial t}v$  and
 $v\cdot\nabla v:=
v_k\frac\partial{\partial x_k}v$. We  investigate on the existence of weak solutions. Here the notion of weak solution (see Definition\,\ref{DWS}) is meant in the sense of Leray's weak solution, but with no change in the arguments it is possible to consider a weak solution in the sense of Caffarelli-Kohn-Nirenberg \cite{CKN}. In this connection we refer the interested reader  to the paper \cite{CGM} for  details.  \begin{defi}\label{DWS}{\sl Let $v_0\in J^2(\OO)$. A field $v:(0,\infty)\times\OO\to \R^3$ is said a weak solution if the field $v$ verifies\begin{itemize}\item[i.]for all $T>0,\,v\in L^\infty(0,T;J^2(\OO))\cap L^2(0,T;J^{1,2}(\OO))$\,,\item[ii.]$\displ\lim_{t\to0}\dm v(t)-v_0\dm_2=0$\,,\item[iii.]for all $t,s\in(0,T)$,  $v$
satisfies the  equation:
$$\displ\intll
st\Big[(v,\varphi_\tau)-(\nabla
v,\nabla
\varphi)+(v\cdot\nabla\varphi,v)\Big]d\tau+(v(s),\varphi
(s))=(v(t),\varphi(t)),$$
$$\mbox{ for all } \varphi\in C^1_0([0,T)\times\OO )\,,\mbox{ with }\n\cdot \varphi=0\,.$$ \end{itemize}}\end{defi} It is well known that any Leray's weak solution enjoys the energy inequality in strong form \cite{Ly,GM,MS}:
\be\label{ESF}\dm v(t)\dm_2^2+\intll st\dm \n v(\tau)\dm_2^2d\tau\leq \dm v(s)\dm_2^2\,,\mbox{ for all }t>s,\mbox{ a.e. in }s>0\mbox{ and for }s=0\,.\ee In the recent paper \cite{CGM} the energy inequality \rf{ESF} becomes a particular case of a suitable energy relation.
In  \cite{CGM}, assuming only $v_0\in J^2(\OO)$,  the existence of a suitable approximating sequence $\{v^m\}$    that, for all $p\in[1,2)$ and $T>0$, enjoys the strong convergence   in $L^p(0,T;J^{1,2}(\OO))\cap L^2(0,T;J^2(\OO))$ to the  weak solution $v$ of \rf{NS} is obtained. This is the first result of strong convergence of the gradient. Its interest is connected with the possibility of deducing an energy equality. This result surely holds  if the convergence is strong in $L^2(0,T;J^{1,2}(\OO))$. Of course  this is not our  case. The idea is to      handle the  strong convergence in $L^p(0,T;J^{1,2}(\OO))$, $p\in[1,2)$, suitably considering  also  the energy relation that holds for the approximation $\{v^m\}$. Unfortunately this it is not enough  for the energy equality. Actually, by our arguments for the field $v$  only   \be\label{EE-I} \dm v(t)\dm_2^2+2\intll st\dm \n v(\tau)\dm_2^2d\tau-\dm v(s)\dm_2^2=-H(t,s), \mbox{a.e. in }t> s>0\mbox{ and for }s=0\,, \ee and  
\be\label{EE-II}  {\displ2\intll st\dm \n v(\tau)\dm_2^2d\tau} = F(t,s)\big({\dm v(s)\dm_2^2-\dm v(t)\dm_2^2}\big), \mbox{a.e. in }t> s>0\mbox{ and for }s=0\, \ee are fulfilled.
 The functions  $H(t,s)$ and $F(t,s)\in(0,1]$ have  the expressions:
\be\label{SEE}\ba{l}H(t,s):=\displ\lim_{\gamma\to0} \lim_{m\to\infty}\gamma\!\intll st\! \frac{\|v^m(\tau)\dm_2^2}{\left(K+\|\nabla v^m(\tau)\|_2^2\right)^{\gamma+1}\hskip-0.3cm}\hskip0.25cm\frac d{d\tau}\|\nabla v^m(\tau)\|_2^2\,d\tau,\mbox{ for  }s>0,\VS H(t,0):=\lim_{s\to0} H(t,s)\ea\ee 
for any arbitrary constant $K>0$, as well as 
 $v$ satisfies relation  \rf{EE-II} with 
\be\label{SEE-I}F(t,s):=\lim_{\gamma\to0}\lim_{m\to\infty}\frac 1{\big(K_1+\dm \n v^m(t_{\gamma,m})\dm_2\big)^\gamma\hskip-0.2cm}\hskip0.2cm\,,\ee 
for any arbitrary constant $K_1\geq0$. Finally,  the following property of inclusion holds:  $\mathcal G_1:=\{t,s\mbox{ such that }\rf{EE-I}\mbox{ is true}\}\subseteq\mathcal G_2:=\{t,s\mbox{ such that }\rf{EE-II}\mbox{ is true}\}$\,. Actually in paper \cite{CGM} these sets of instants are the ones for which the strong convergence of the sequence $\{v^m\}$ holds in $J^{1,2}(\OO)$ and in $J^2(\OO)$, respectively. In the proof of formulas \rf{EE-III} and \rf{EE-IV} we also use the set $\mathcal G_1$ as set of instants for which the quoted strong convergence holds.  \par In this note we give further partial contributions for an energy equality that are of some interest. Firstly we provide a generalization of properties \rf{EE-I} and \rf{EE-II} respectively. Actually we are able to prove
\be\label{EE-III}G(t,s)=  \dm v(s)\dm_2^2-\dm v(t)\dm_2^2-2\intll st\dm \n v(\tau)\dm_2^2d\tau\, , \mbox{a.e. in }t\geq s>0\mbox{ and for }s=0\,, \ee  and  \be\label{EE-V}  {\displ2\intll st\dm \n v(\tau)\dm_2^2d\tau} = P(t,s)\big({\dm v(s)\dm_2^2-\dm v(t)\dm_2^2}\big), \mbox{a.e. in }t\geq s>0\mbox{ and for }s=0\,. \ee 
   The functions  $G(t,s)$ and $P(t,s)$ have  the expressions:
\be\label{SEE-III}\ba{l}G(t,s)\!:=-\displ\lim_{\beta\to0}\lim_{m\to\infty}\beta\!\intll st\! \dm v^m(\tau)\dm_2^2\, g^{\beta-1}(\dm \n v^m(\tau)\dm_2^2)g'(\dm \n v^m(\tau)\dm_2^2)\frac d{d\tau}\|\nabla v^m(\tau)\|_2^2\,d\tau,\mbox{\,for\,}s\!>\!0,\VS G(t,0)\!:=\lim_{s\to0}G(t,s),\ea\ee 
for any arbitrary function $g(\rho)>0$ which is Liptchiz continuous on $[0,\infty)$ and such that $\displ\lim_{\rho\to\infty}g(\rho)=0$, as well as 
 $v$ satisfies relation  \rf{EE-V} with alternative expressions  for $P(t,s)$ \be\label{SEE-V}\lim_{R\to\infty}p_R^\alpha(h(\alpha,R))\in(0,1]\,,\;\lim_{\alpha\to0^+}\lim_{R\to\infty}p_R^\alpha(h(\alpha,R))\in(0,1]\,,\;\lim_{\alpha\to\infty}\lim_{R\to\infty}p_R^\alpha(h(\alpha,R))\in(0,1]\,,\ee where the exponent $\alpha$ is  real, $p_R:[0,\infty)\to [0,1]$ is a positive continuous cut off function with $p_R(\rho)=1$ for $\rho\in[0,R]$, $p_R(\rho)=0$ for $\rho\geq2R$, and $p_R(\rho)\in (0,1)$ for $\rho\in (R, 2R)$, and finally $h(\alpha,R)\in [0,2R)$. Of course by proving that $G(t,s)\leq0$ or $P(t,s)=1$ one deduces the energy equality. No contraindication is known for these claims.  Moreover we give the following\be\label{EE-IV} E(t,s)  =  {\dm v(s)\dm_2^2-\dm v(t)\dm_2^2}-{\displ2\intll st\dm \n v(\tau)\dm_2^2d\tau} \,, \mbox{a.e. in }t\geq s>0\mbox{ and for }s=0\,, \ee with
\be\label{SEE-IV}\hskip-0.5cm\ba{l}E(t,s):=\displ\lim_{R\to\infty}E(t,s,R)\VS E(t,s,R)\!:=\!\!\displ \lim_{m\to\infty}\!\left[\mbox{${\underset{h=1}{\overset {p_{22}(m,R)}\sum}}$}\!2\Big[\dm  v^m(t_h^{22})\dm_2^2-\dm v^m(s_h^{22})\dm_2^2\Big] \!+\!\!\mbox{${\underset{h=1}{\overset {p_{21}(m,R)}\sum}}$}\Big[\dm  v^m(t_{h}^{21})\dm_2^2-2\dm v^m(s_h^{21})\dm_2^2\Big]\right.\VS\left.\hskip3cm \!\!\!+\!\mbox{${\underset{h=1}{\overset {p_{11}     (m,R)}\sum}}$}\Big[\dm  v^m(t_h^{11})\dm_2^2-\dm v^m(s_h^{11})\dm_2^2\Big]\! +\!\mbox{${\underset{h=1}{\overset {p_{12}(m,R)}\sum}}$}\Big[2\dm  v^m(t_h^{12})\dm_2^2-\dm v^m(s_h^{12})\dm_2^2\Big]\right]\VS E(t,0):=\lim_{s\to0}E(t,s)\,,\ea\ee where, for all $R>0$\,, $t_h^{ij}>s_h^{ij}$ and \be\label{MIS} \mbox{${\underset{h=1}{\overset {p_{22}(m,R)}\sum}}$}\!(t_h^{22}\!-s_h^{22})+\!
\mbox{${\underset{h=1}{\overset {p_{12}(m,R)}\sum}}$}\!(t_h^{12}\!-s_h^{12})
+\!\mbox{${\underset{h=1}{\overset {p_{11}(m,R)}\sum}}$}\!(t_h^{11}\!-s_h^{11})+\!\mbox{${\underset{h=1}{\overset {p_{21}(m,R)}\sum}}$}\!(t_h^{21}\!-s_h^{21})<\frac {\dm v_0\dm_2^2}{2R}\,.\ee 
Since $$\dm v^m(t_h^{ij})\dm_2-\dm v^m(s_h^{ij})\dm_2<0,\mbox{ for all }{i,j, h}\,,$$ we get that   on the right-hand side of \rf{SEE-IV} only the last sum can be positive, the other sums surely are negative. In addition, by  construction   the equality \be\label{SIGN}p_{21}(m,R)=p_{12}(m,R)\ee holds. Hence we get \be\label{SUM}E(t,s,R)=\lim_{m\to\infty}\Big[B(m,R)+\dm v^m(t_1^{12})\dm_2^2-\dm v^m(s_{p_{21}}^{21})\dm_2^2\Big]\,,\mbox{ with }B(m,R)<0\,,\mbox{ for all }R>0\,.\ee  This is of some interest. Actually if $E(t,s)\leq0$, then we trivially have that the energy equality holds.  \par In the linear case or for 2D-Navier-Stokes equations we get $G(t,s)=0$, $P(t,s)=1$ and $E(t,s)=0$. The first two claims are already proved in \cite{CGM}.  In the case of $E(t,s)$, for $R$ sufficiently large, we are able to prove that $p_{ij}(m,R)=0$ for all $i,j=1,2$. \par All the above considerations are condensed in the following statement which is the chief result of the note:\begin{tho}\label{CT}{\sl For all $v_0\in J^2(\OO)$ there exists a  weak solution $v$ to problem \rf{NS} which is a weakly continuous function of $t\geq0$. The field $v$ is the weak* limit in $L^\infty(0,T;J^2(\OO))$ and weak limit in $L^2(0,T; J^{1,2}(\Omega))$ of a sequence $\{(v^m,\pi_{v^m})\}$ of solutions to \eqref{MNS}.
The sequence $\{v^m\}$ strongly converges   to $v$ in $L^2(0,T;L^2(\OO))\cap L^p(0,T;W^{1,2}(\OO))$ for all $p\in [1,2)$. Moreover the limit $v$ satisfies the energy relations \rf{EE-III}-\rf{EE-IV} with functions $G(t,s)$, $P(t,s)$  and $E(t,s)$ given by \rf{SEE-III}-\rf{SEE-IV} respectively. In particular  a.e. in $s> 0$ and in $s=0$,  $v$ is continuous  on the right  in the $L^2$-norm. Finally, in the case of equality \rf{EE-IV} estimates \rf{MIS}-\rf{SUM} hold. }\end{tho} \begin{rem}{\rm By virtue of \rf{EE-III}-\rf{EE-IV}  the energy inequality \rf{ESF} is a particular case again. We stress that if the energy equality does not  hold, then the gap with the equality is governed by the terms $G(t,s),P(t,s)$ and $E(t,s)$. Since for the definition of $G$ we have arbitrary auxiliary     Lipschitz continuous functions $g(\rho)$ and for the definition of $P(t,s)$ arbitrary exponents $\alpha$ and its limits to $0$ or $\infty$ respectively,  we have a wide set of expressions related to the same value of the gap. In the case of the term $E(t,s)$, by means of a possible series, the gap is expressed by means of the $L^2$-norm of the approximations. \par We remark that in the case of a bounded domain $\OO$ the results also hold assuming  the Galerkin approximation as sequence, as proposed in \cite{Pr} and in \cite{Hy}.\par Finally, if we consider the validity of the energy equality for the weak limit, for example by assuming, for all $\vep>0$, $v\in L^4(\vep,T;L^4(\OO))$, as made in \cite{M-I}, which is compatible with the assumption $v_0\in J^2(\OO)$, then we get $G(t,s)=E(t,s)=0$ and $P(t,s)=1$. It is of some interest to understand if in turn these new relations can imply an information on the partial  regularity or on the  singularity of the solution $v$, limit of the sequence $\{v^m\}$ quoted in Theorem\,\ref{CT}.  }\end{rem} 
\section{Notations and some preliminary result.}
The symbol $\mathscr C_0(\OO)$ indicates the set of function $C_0^\infty(\OO)$ and divergence free. In this paper we denote by $J^2(\OO)$ and $J^{1,2}(\OO)$ the completion spaces of $\mathscr C_0(\OO)$ with respect to the metric of $L^2(\OO)$ and $W^{1,2}(\OO)$ respectively. By $L^p(0,T;X)$ we mean the space of functions 
from $(0,T)$ into $X$ such that $\intll0T\dm u(\tau)\dm_X^pd\tau<\infty$. The following lemma is not known to the authors, hence the proof is proposed. \begin{lemma}\label{LEBT}{\sl Let $\{h^m(t)\}\subset L^1((0,T))$ be a bounded sequence. Assume that $h^m(t)\to h(t)$ a.e. in $t\in(0,T)$ and $h\in L^1((0,T))$. Let $p(s)\in C([0,\infty))$ with $p(s)>0$ and $\displ\lim_{s\to\infty}p(s)=0$. Then, there exists a function $h_0(t)\equiv h(t)$ a.e. in $t\in(0,T)$, such that, for all $\alpha>0$, we get
\be\label{LP}\lim_m\intll0T\big|{h^m(t)}{p^\alpha(|h^m(t)|)}-{h(t)}{p^\alpha(|h_0(t)|)} \big|dt=0\,.\ee Moreover,  the following limit property also holds:\be\label{LA}\lim_{\alpha\to0}\lim_m\intll0Th^m(t)p^\alpha(|h^m(t)|)dt=\intll0T h(t)dt\,.\ee}\end{lemma}
\bp Since $h\in L^1((0,T))$, the estimate $|h(t)|<\infty$ holds for all $t\in (0,T)-A$ with $meas(A)=0$\,.   We set $h_0(t)=h(t)$ for  $t\in (0,T)-A$ and $h_0(t)=0$ for $t\in A$.           
We  write the left-hand side of \rf{LP} as follows:\be\label{SLM-I}\ba{ll}\displ\intll0T\hskip0.5cm&\displ\hskip-0.7cm\big|{h^m(t)}{p^\alpha(|h^m(t)|)}-{h(t)}{p^\alpha(|h_0(t)|)} \big|dt\\&\displ =\intll0T\big|\big[{h^m(t)-h(t)}\big]{p^\alpha(|h^m(t)|)}+{h(t)}\big[{p^\alpha(|h^m(t)|)}-{p^\alpha(|h_0(t)|)}\big]\big|dt\\&\displ\leq \intll0T{|h^m(t)-h(t)|}{p^\alpha (|h^m(t)|)}dt+\intll0T{|h(t)|}\big|p^\alpha(|h^m(t)|)-p^\alpha(|h_0(t)|)\big|dt \,.\ea\ee
Since $\big|p^\alpha(|h^m(t)|)-p^\alpha(|h_0(t)|)\big|\leq 2\dm p\dm_\infty^\alpha\,,$ for all $t\in(0,T)$, and $h^m(t)\to h(t)$  a.e. in $t\in(0,T)$, applying the Lebesgue dominated convergence  theorem, we arrive at \be\label{IP}
 \lim_m \intll0T{|h(t)|}\big|p^\alpha(|h^m(t)|)-p^\alpha(|h_0(t)|)\big|dt =0\,.\ee Concerning the first integral on the right-hand side of \rf{SLM-I}, we arrange  the following artifice:$$\ba{ll}\displ\intll0T{\!|h^m(t)\!-h(t)|}{p^\alpha(|h^m(t)|)}dt\hskip-0.2cm&\displ\!=\!\!\intl{I^m_M}{\!|h^m(t)\!-h(t)|}{p^\alpha(|h^m(t)|)}dt+\!\hskip-0.7cm\intl{\hskip0.2cm(0,T)-I^m_M}\hskip-0.6cm{\!|h^m(t)\!-h(t)|}{p^\alpha(|h^m(t)|)}dt\VSE=:J_1^m+J_2^m\,,\ea$$ where, for some fixed $M>0$,  $I^m_M:=\{t:|h^m(t)|\leq M\}$\,.
By the definition of $I^m_M$, we easily obtain $$J_2^m\leq p^\alpha_M\intll0T(|h^m(t)|+|h(t)|)dt\leq  c{p^\alpha_M}\,,\mbox{ for all }m\in\N\,,$$ where $c$ is independent of $M$  and  $\displ p_M:=\sup_{(M,\infty)} p(s)$\,. Concerning $J^m_1$, we introduce the  characteristic function of the set $I^m_M$, so that$$\ba{c}H^m(t):={|h^m(t)-h(t)|}{p^\alpha(|h^m(t)|)}\chi_{I^m_M}(t)\to 0 \mbox{ a.e. in }t\in(0,T)\,,\VS\mbox{ and,  a.e. in }t\in(0,T)\,,\;H^m(t)\leq (M+|h(t)|)\dm p\dm_\infty^\alpha\mbox{ for all }m\in\N \,.\ea$$ So that applying the Lebesgue dominated convergence  theorem, we deduce that $\displ\lim_mJ^m_1\!=0\,.$ The above limit property of $J_1^m$ and estimate on $J_2^m$,  via \rf{SLM-I}-\rf{IP}, furnish
$$\lim_m\intll0T\big|{h^m(t)}{p^\alpha(|h^m(t)|)}-{h(t)}{p^\alpha(|h_0(t)|)} \big|dt\leq  c{p^\alpha_M},\mbox{ for all }M>0\,.$$ Since $p\to0$ for $s\to\infty$, letting $M\to\infty$ we prove \rf{LP}\,. By applying the Lebesgue dominated convergence  theorem we   immediately get \rf{LA}\,.
\ep
\begin{coro}\label{LEBTT}{\sl Let $\{h^m(t)\}\subset L^1((0,T))$ be a  sequence with $A:=\displ\sup_m\dm h^m\dm_1<\infty$\,. Assume that $h^m(t)\to h(t)$ a.e. in $t\in(0,T)$ and $h\in L^1((0,T))$. Let $p_R(s)\in C([0,\infty))$ with $p_R(s)=1$ for $s\leq R$, $p_R(s)=0$ for $s\geq 2R$ and $p_R(s)\in (0,1)$ for $s\in (R, 2R)$. Then, there exists a function $h_0(t)\equiv h(t)$ a.e. in $t\in(0,T)$, such that, for all $R>(2T)^{-1}A$, we get
\be\label{LC}\lim_m\intll0T\big|{h^m(t)}{p_R(|h^m(t)|)}-{h(t)}p_R(|h_0(t)|) \big|dt=0\,.\ee Moreover,  the following limit also holds:\be\label{LCA}\lim_{R\to\infty}\lim_m\intll0Th^m(t)p_R(|h^m(t)|)dt=\intll0T h(t)dt\,.\ee}\end{coro}
\bp The first claim of the corollary is a particular case of the previous lemma.  Hence we get \be\label{SLM-V}\ba{ll}\displ\intll0T\hskip0.5cm&\displ\hskip-0.7cm\big|{h^m(t)}{p_R(|h^m(t)|)}-{h(t)}{p_R(|h_0(t)|)} \big|dt\\&\displ =\intll0T\big|\big[{h^m(t)-h(t)}\big]{p_R(|h^m(t)|)}+{h(t)}\big[{p_R(|h^m(t)|)}-{p_R(|h_0(t)|)}\big]\big|dt\\&\displ\leq \intll0T{|h^m(t)-h(t)|}{p_R (|h^m(t)|)}dt+\intll0T{|h(t)|}\big|p_R(|h^m(t)|)-p_R(|h_0(t)|)\big|dt \,,\ea\ee
and  \be\label{LC-I}
 \lim_m \intll0T{|h(t)|}\big|p_R(|h^m(t)|)-p_R(|h_0(t)|)\big|dt =0\,.\ee Concerning the first integral on the right-hand side of \rf{SLM-V}, we get $$\displ\intll0T{\!|h^m(t)\!-h(t)|}{p_R(|h^m(t)|)}dt\displ\!=\!\!\intl{T_R}{\!|h^m(t)\!-h(t)|}{p_R(|h^m(t)|)}dt=:J^m_R\,,$$ where we set $T_R^m:={\{t:0\leq h^m(t)<2R\}}$\,\footnote{\;We stress that, for all $m\in\N$, $meas(T_R^m)\ne0$. Actually, if for $2TR>A$ we guess that $meas(T_R^m)=0$, then $$A<\intll0T2R\leq\intll0T h^m(t)dt\leq A\,.$$}. 
Since  $$\ba{c}H^m(t):={|h^m(t)-h(t)|}{p_R(|h^m(t)|)}\to 0 \mbox{ a.e. in }t\in T_R^m\,,\VS H^m(t)\leq 2R\,,\;\mbox{ for all }m\in\N\mbox{ and }t\in T_R^m \,,\ea$$  applying the Lebesgue dominated convergence theorem, we deduce that, for all $R>0$, letting $m\to\infty$, $J^m_R\to 0$.  The above limit property  via \rf{SLM-V}-\rf{LC-I} furnishes
$$\lim_m\intll0T\big|{h^m(t)}{p_R(|h^m(t)|)}-{h(t)}{p_R(|h_0(t)|)} \big|dt=0\,,$$\mbox{ hence } $$\lim_m\intll0T{h^m(t)}{p_R(|h^m(t)|)}=\intll0T{h(t)}{p_R(|h_0(t)|)}dt dt\,.$$  The first limit proves \rf{LC}. By applying  Lebesgue's theorem to the second  identity, letting $R\!\to\!\infty$, we   immediately get \rf{LCA}\,.
\ep We consider  a mollified Navier-Stokes initial boundary value problem:\be\label{MNS}\ba{l}v_t^m+J_m[v^m]\cdot
\nabla v^m+\nabla\pi_{v^m}=\Delta
v^m,\;\nabla\cdot
v^m=0,\mbox{ in }(0,T)\times\OO,\VS v^m=0\mbox{ on }(0,T)\times\po,\hskip0.12cm
v^m(0,x)=v_0^m(x)\mbox{ on
}\{0\}\times\OO,\ea\ee where $J_m[\cdot]$ is a mollifier and  $\{v_0^m\}\subset J^{1,2}(\OO)$ converges to $v_0$ in $J^2(\OO)$. In \cite{CGM} is proved 
\begin{lemma}\label{LJAM}{\sl For all $v_0\in J^2(\OO)$ there exists a sequence  of solutions $\{(v^m,\pi_{v^m})\}$ such that, for all $m\in\N$ and $T>0$, $v^m\in C([0,T);J^{1,2}(\OO))\cap L^2(0,T;W^{2,2}(\OO))$.}\end{lemma}\begin{lemma}\label{EX}{\sl The sequence  $\{v^m\}$ of kinetic fields of Lemma\,\ref{LJAM}  admits  limit $v$ which  is weak$^*$ limit in $L^\infty(0,T;J^2(\OO))$, weak limit in  $L^2(0,T;J^{1,2}(\OO))$ and, for all $p\in[1,2)$, strong limit in $L^p(0,T;J^{1,2})\cap L^2(0,T;J^2(\OO))$.  Moreover,   $v$ is a weak solution to problem \rf{NS} with $(v(t),\psi)$ continuous function of $t$, for all $\psi\in J^2(\OO)$.}\end{lemma}\bp Apart from the strong convergence properties,  the result  is well known. The properties of strong convergence are proved in \cite{CGM}.\ep 
\begin{lemma}\label{CCM}{\sl Let $(\{v^m,\pi_{v_m})\}$ be the sequence furnished by Lemma\,\ref{LJAM}. Then, for all $t\geq0$ and $m\in\N$, $\dm v^m\dm_{1,2}\ne0$ holds.}\end{lemma} \bp See Lemma\,10 in \cite{CGM}.\ep\begin{lemma}\label{CC}{\sl Let $v$ be a weak solution. Assume that  $ \dm v(s)\dm_2\ne 0$ and $v$ is right continuous in the $L^2$-norm. Then there exists $\delta>0$ such that $\dm v(t)\dm_2\ne0$ for all  $t\in[s,s+\delta)$.}\end{lemma}\bp The proof is given in \cite{CGM}. However since it is very short,   for the sake of the completeness we reproduce it. Assume that for all $\delta>0$ we have a $t\in(s,s+\delta)$ such that $\dm v(t)\dm_2=0$. Then  we can select  a sequence $\{t_p\}$ converging to $s$ such that $\dm v(t_p)\dm=0$. Then,  by virtue of the right-$L^2$-continuity in $s$, we get $\lim_{t_p\to s^+}\dm v(t_p)-v(s)\dm_2= 0$, which is a contradiction with $\dm v(s)\dm_2\ne 0$.\ep\section{Proof of Theorem\,\ref{CT}} Taking into account the results of Lemma\,\ref{LJAM} and Lemma\,\ref{EX}, our proof is reduced  to prove the energy inequalities \rf{EE-III}-\rf{EE-V}. In the sequel by $\mathcal G_1$ we mean the set $$\{t\geq0:\dm v^m(t)-v(t)\dm_{1,2}\to0\mbox{ as }m\to\infty\}\,.$$ Thanks to Lemma\,\ref{EX} we have $meas(\R^+-\mathcal G_1)=0$. In the sequel by $\mathcal G_2$ we mean the set $$\{t\geq0:\dm v^m(t)-v(t)\dm_{2}\to0\mbox{ as }m\to\infty\}\,.$$ We have $meas(\R^+-\mathcal G_2)=0$, and   {\it a priori }$\mathcal G_1\subseteq\mathcal G_2$ holds.
\par {\it Energy equality \rf{EE-III}.}\newline The starting point is the energy differential equation of $(v^m,\pi_{v^m})$ related to the solution of the  problem \rf{MNS}, that is 
\be\label{P-O} \frac d{dt}\dm v^m(t)\dm_2^2+2\dm\n v^m(t)\dm_2^2=0\,,\;t\geq0\,. \ee
We multiply equation \rf{P-O} by $g^\beta(\dm\n v^m(t)\dm_2^2)$. Here $g(\rho)>0$ is a Lipschitz continuous function of $\rho\geq0$ such that $\displ\lim_{\rho\to\infty}g(\rho)=0$ and $\beta>0$ is an exponent. We consider $s,t\in \mathcal G_1$ with $0<s<t$. Integrating by parts  on the interval $(s,t)$, we get \be\label{P-I}\hskip-0.3cm\ba{l}\displ -\beta\intll st\dm v^m(\tau)\dm_2^2g^{\beta-1}(\dm \n v^m(\tau)\dm_2^2)g'(\dm \n v^m(\tau)\dm_2^2)\frac d{dt}\dm \n v^m(\tau)\dm_2^2d\tau\\\displ\hskip1.1cm=\dm v^m(s)\dm_2^2g^\beta(\dm\n v^m(s)\dm_2^2)-\dm v^m(t)\dm_2^2g^\beta(\dm \n v^m(t)\dm_2^2)-2\!\intll st\!g^\beta(\dm\n v^m(\tau)\dm_2^2)\dm \n v^m(\tau)\dm_2^2d\tau .\ea\ee Recalling the definition of $\mathcal G_1$, by virtue of Lemma\,\ref{LEBT}, letting $m\to\infty$ and subsequently letting $\beta\to0$, we arrive at \rf{EE-III} with $G(t,s)$ given by \rf{SEE-III}.\chiu\par{\it Energy equality \rf{EE-V}}.\newline  
We consider $s,t\in \mathcal G_2$ with $0\leq s<t$. We introduce  a cutoff function $p_R(\rho)$ which is Lipschitz continuous and such that $p_R(\rho)=1$ for $\rho\in[0,R]$, $p_R(\rho)=0$ for $\rho\geq2R$ and $p_R(\rho)\in (0,1)$ for $\rho\in (R, 2R)$.  We set \be\label{RA}A=\intll 0t\dm \n v(\tau)\dm_2^2d\tau, \mbox{ and }R_0>(2t)^{-1}A\mbox{ and such that }meas(\{\tau\in(0,t):\dm \n v(\tau)\dm_2^2<2R_0\})\ne0\,.\ee The claim \rf{RA} is consistent  by virtue of Lemma\,\ref{CC}. For $\alpha>0$, we get \be\label{CFE}p_R^\alpha(\dm\n v^m(t)\dm_2^2)\mbox{$\frac 12\frac d{dt}$}\dm   v^m(t)\dm_2^2+p_R^\alpha(\dm \n v^m(t)\dm_2^2)\dm \n v^m(t)\dm_2^2=0\,.\ee Integrating on $(0,t)$, considering $R>R_0$ and applying the mean value theorem, for some $t_{\alpha,R}^m\in(0,t)$, we get
$$p_R^\alpha(\dm\n v^m(t^m_{\alpha,R})\dm_2^2)\mbox{$\frac 12$}\big[\dm   v^m(t)\dm_2^2-\dm v^m(0)\dm_2^2\big]=-\intll0tp_R^\alpha(\dm \n v^m(\tau)\dm_2^2)\dm \n v^m(\tau)\dm_2^2d\tau\,.$$ By virtue of Lemma\,\ref{CCM}, for all  $m\in\!\N$ the right-hand side is different from to zero, we deduce that $\dm\n v^m(t^m_{\alpha,R})\dm_2^2\!\in (0,2R)$ for all $m\in\N$. Hence there exists an extract which admits a limit $h(\alpha,R)\in [0,2R]$. We label again by $m$  the extract sequence.  Letting $m\to\infty$, firstly  we obtain$$\lim_m p_R^\alpha(\dm\n v^m(t^m_{\alpha,R})\dm_2^2)=\frac{2\intll0t p_R^\alpha(\dm \n v(\tau)\dm_2^2)\dm \n v(\tau)\dm_2^2d\tau}{\dm   v(0)\dm_2^2-\dm v(t)\dm_2^2}\ne0\,,$$ as a consequence of \rf{RA}.  Then,  since, for all $\alpha>0$ and $R>R_0$, $p^\alpha_R$ is continuous, we get $\displ \lim_m p^\alpha_R(\dm\n v^m(t^m_{\alpha,R})\dm_2^2)=p^\alpha_R(h(\alpha,R))$. So that  letting $R\to\infty$, we get
$$ \lim_Rp^\alpha_R(h(\alpha,R))= \frac{2\intll0t \dm \n v(\tau)\dm_2^2d\tau}{\dm   v(0)\dm_2^2-\dm v(t)\dm_2^2}\ne0\,,$$ which proves \rf{EE-V} with $P(t,s)$ given by \rf{SEE-V}$_1$. The alternative expressions of $P(t,s)$  are an immediate consequence of the fact that the right-hand side is independent of $\alpha$\,. \chiu\par{\it Energy equality \rf{EE-IV}}.\newline We consider $s,t \in \mathcal G_1$ with $0<s<t$. We consider equation \rf{CFE} with $p_R(\rho):=\frac{2R-\rho}R\,,$ where $R>\max\{\dm \n v(s)\dm_2^2,\dm \n v(t)\dm_2^2\}$\,. Since $s,t\in\mathcal G_1$ we have strong convergence in $J^{1,2}(\OO)$ of $\{v^m(s,x)\}$ to  $v(s,x)$ and of $\{v^m(t,x)\}$ to $v(t,x)$, respectively. Hence, without losing the generality,  we also claim that $R>\max\{\dm \n v^m(s)\dm_2^2,\dm \n v^m(t)\dm_2^2\}$ for all $m\in\N$\,.  Integrating on $(s,t)$, after an integration by parts we get
\be\label{CR}-\intll st \dm v^m(\tau)\dm_2^2\frac d{d\tau}p_R(\dm \n v^m(\tau)\dm_2^2)d\tau={\dm v^m(s)\dm_2^2-\dm v^m(t)\dm_2^2}-{\displ2\intll stp_R(\dm\n v^m(\tau)\dm_2^2)\dm \n v^m(\tau)\dm_2^2d\tau} \,.\ee
We set   $$I^m_R:=\{\tau\in(s,t):R<\dm \n v^m(\tau)\dm_2^2<2R\}\,.$$ The set $I_R^m$, if not empty, is at most a countable union of disjoint open intervals. Actually, let $J$ be the interior set of $(s,t)-I^m_R$. Since $\partial I_R^m$ is at most countable, we get $$\intll st\! \dm v^m(\tau)\dm_2^2 \frac d{d\tau} p_R(\dm \n v^m(\tau)\dm_2^2)d\tau =\!\!\intl {\hskip-0.2cmI_R^m} \! \dm v^m(\tau)\dm_2^2 \frac d{d\tau} p_R(\dm \n v^m(\tau)\dm_2^2)d\tau+\intl J\!\dm v^m(\tau)\dm_2^2 \frac d{d\tau} p_R(\dm \n v^m(\tau)\dm_2^2)d\tau\,.$$ Let us fix $\ov \tau\in J$ and suppose that $\dm \n v^m(\ov \tau)\dm_2^2\leq R$. Since $J$ is open, it contains an interval centered at $\ov\tau$ where, by continuity, $\dm \n v^m(\tau)\dm_2^2\leq R$ holds. Hence we get $\frac d{d\tau}p_R(\dm v^m(\tau)\dm_2^2)=0$\,. The same argument works in the case of $\dm \n v^m(\ov\tau)\dm_2^2\geq2R$\,. Thus, we finally get 
$$\intll st \dm v^m(\tau)\dm_2^2\frac d{d\tau}p_R(\dm \n v^m(\tau)\dm_2^2)d\tau= \intl {I_R^m}\dm v^m(\tau)\dm_2^2 \frac d{d\tau} p_R(\dm \n v^m(\tau)\dm_2^2)d\tau\,,$$ which coupled with \rf{CR} furnish    $$\displ\intl{I^m_R}\frac{\dm v^m(\tau)\dm_2^2}R\frac d{d\tau}\dm \n v^m(\tau)\dm_2^2d\tau ={\dm v^m(s)\dm_2^2-\dm v^m(t)\dm_2^2}-{\displ2\intll stp_R(\dm\n v^m(\tau)\dm_2^2)\dm \n v^m(\tau)\dm_2^2d\tau} \,.$$  We denote by $p_{ij}(m,R),$ $i,j=1,2$, a positive or null integer. If we have $i\ne j$, then $p_{ij}(m,R)\in\N$ (see item 3. below). If we have $i=j$, then $p_{ii}(m,R)$ can be $\infty$\,. Actually, we distinguish the intervals of the union in the following way:
\begin{itemize}
\item  $h\!\in\{1,\dots,p_{11}(m,R)\}$, where {\it a priori} $p_{11}(m,R)\in\N\cup\{\infty\}$, $(s^{11}_h,t^{11}_h)$ is the interval such that $\dm \n v^m(\tau)\dm_2^2\in(R,2R)$\,, for all $\tau\in(s_h^{11},t_h^{11})$\,, with $\dm\n v^m(s^{11}_h)\dm_2^2=\dm\n v^m(t^{11}_h)\dm_2^2=R$\,,
\item\,$h\!\in\{1,\dots,p_{12}(m,R)\}$, $(s^{12}_h,t^{12}_h)$ is the interval such that $\dm\n v^m(\tau)\dm_2^2\in(R,2R)$\,, for all $\tau\in(s_h^{12},t_h^{12})$\,, with $\dm \n v^m(s^{12}_h)\dm_2^2=R$  and $\dm\n v^m(t^{12}_h)\dm_2^2=2R$\,,
   \item  $h\!\in\{1,\dots,p_{21}(m,R)\}$\,, $(s^{21}_h,t^{21}_h)$ is the interval such that $\dm\n v^m(\tau)\dm_2^2\in(R,2R)$\,, for all $\tau\in(s_h^{21},t_h^{21})$\,, with $\dm \n v^m(s^{21}_h)\dm_2^2=2R$ and $\dm\n v^m(s^{21}_h)\dm_2^2=R$\,, 
\item  $h\!\in\{1,\dots,p_{22}(m,R)\}$, where {\it a priori} $p_{22}(m,R)\in\N\cup\{\infty\}$,\,$(s^{22}_h,t^{22}_h)$ is the interval such that $\dm \n v^m(\tau)\dm_2^2\in(R,2R)$\,, for all $\tau\in(s_h^{22},t_h^{22})\,,$ with $\dm\n v^m(s^{22}_h)\dm_2^2=\dm \n v^m(t^{22}_h)\dm_2^2=2R$\,.
\end{itemize}
Since $\dm \n v^m(s)\dm_2^2<R$ and $\dm\n v^m(t)\dm_2^2<R$,  we deduce that $s<\min s^{ij}_h$ and $t>\max t^{ij}_h\,.$ Moreover
\begin{enumerate}\item $p_{ij}(m,R)=0$, for all $i,j\in \{1,2\}$, for all $m\in\N$ if, and only if, $\dm \n v^m(\tau)\dm_2^2<R$ for all $m\in\N$\,.  \item $p_{12}(m,R)=p_{21}(m,R)=p_{22}(m,R)=0$ if, and only if, $\dm \n v^m(\tau)\dm_2^2< 2R$ for all $\tau\in(s,t)$\,.\item  $p_{12}(m,R)\ne0$ if, and only if, $p_{21}(m,R)\ne0$\,, and in both  cases $p_{ij}(m,R)\in \N$\,. Actually,  we have that $\{(s_h^{12},t_h^{12})\}$ is a sequence, the sequence collapses into a point $\ov t$ and as a consequence we get $\displ R=\lim_{s_h^{12}\to\ov t}p_R(\dm \n v^m(s_h^{12})\dm_2^2)=p_R(\dm\n v^m(\ov t)\dm_2^2)=\lim_{t_h^{12}\to\ov t}p_R(\dm \n v^m(t_h^{12})\dm_2^2)=2R $ which is a contradiction. Analogous argument works  in the case of ``21''.\item We have  $p_{12}(m,R)=p_{21}(m,R)$. Actually, for fixed $h\leq p_{12}(m,R)$,  let us consider the set $A:=\{k\in\N:t^{12}_h\leq s^{21}_k\}$. Since $\dm \n v^m(t)\dm_2^2<R$, by continuity, the set $A$ is not empty. Let be $\varphi(h):=\min(A)$. Then, by construction, $\dm\n v^m(\tau)\dm_2^2>R$ for all $\tau$ such that $\tau\in (t^{12}_h,t^{21}_{\varphi(h)})$. It follows that $$s_{\varphi( h)}^{21}<t^{21}_{\varphi( h)}<s_{h+1}^{12}<t^{12}_{h+1}<s^{21}_{\varphi({h+1})}\,.$$ Hence function $\varphi$ is injective.  Conversely, using the fact that  $\dm\n v^m(s)\dm_2^2<R$,  we can conclude that $p_{12}(m,R)=p_{21}(m,R)$\,. \end{enumerate}
  Since $\,2\intll st\dm  \n v^m(\tau)\dm_2^2d\tau\leq \dm v_0\dm_2^2\,,$ then, for all $m\in\N$ and $R$, we get $meas(I_R^m)<\frac{\dm v_0\dm_2^2}{2R}$ which proves \rf{MIS}. By the energy inequality the same property holds for the weak solution $v$, that is  \be\label{MIR} I_R:=\{\tau\in (s,t): R\leq \dm \n v(\tau)\dm_2^2\leq 2R\}\,,\;\mbox{then }meas(I_R)  <\frac{\dm v_0\dm_2^2}{2R} \,.\ee       Taking into account the previous considerations, integrating by parts on the intervals defining $I^m_R$, we get \be\label{MEE}E^m= -\frac2R\intl{I_R^m}\dm\n v^m(\tau)\dm_2^4d\tau +{\dm v^m(s)\dm_2^2-\dm v^m(t)\dm_2^2}-{\displ2\intll stp_R(\dm\n v^m(t)\dm_2^2)\dm \n v^m(\tau)\dm_2^2d\tau} \,,\ee 
where we considered that $\frac d{dt}\dm v^m(t)\dm_2^2=-2\dm \n v^m(t)\dm_2^2$, and  we  set
\be\label{DIFF}\ba{l}E^m:=\mbox{${\underset{h=1}{\overset {p_{22}(m,R)}\sum}}$}\!2\Big[\dm  v^m(t_h^{22})\dm_2^2-\dm v^m(s_h^{22})\dm_2^2\Big] \!+\!\!\mbox{${\underset{h=1}{\overset {p_{21}(m,R)}\sum}}$}\Big[\dm  v^m(t_{h}^{21})\dm_2^2-2\dm v^m(s_h^{21})\dm_2^2\Big]\VS\hskip3cm \!\!\!+\!\mbox{${\underset{h=1}{\overset {p_{11}(m,R)}\sum}}$}\Big[\dm  v^m(t_h^{11})\dm_2^2-\dm v^m(s_h^{11})\dm_2^2\Big]\! +\!\mbox{${\underset{h=1}{\overset {p_{12}(m,R)}\sum}}$}\Big[2\dm  v^m(t_h^{12})\dm_2^2-\dm v^m(s_h^{12})\dm_2^2\Big]\,.\ea\ee In the expression of $E^m$ we have   taken into account  that the value of $\dm\n v^m(\tau)\dm_2^2$ in the end points of the intervals is alternatively $R$ or $2R$ and for all  the  integral term there is the moltiplicative  factor $\frac 1R$\,. We discuss the first integral on the right-hand side of \rf{MEE}. By the definition of $I^m_R$ we deduce that $\frac {\dm \n v^m(t)\dm_2^2}R\leq 2$ for all $\tau\in I^m_R$ and for all $m\in\N$\,.
Hence, employing the energy relation for $v^m$ and the Lebesgue theorem, we deduce that for all $R$ $$\lim_{m}\frac 1R\intl{I^m_R}\dm\n v^m(\tau)\dm_2^4d\tau\leq2 \lim_m\intl{I^m_R}\dm \n v^m(\tau)\dm_2^2d\tau\leq 2\intl{I_R}\dm \n v(\tau)\dm_2^2d\tau\,.$$ Since for all $T>0$ we have $\dm\n v(\tau)\dm_2^2\in L^1(0,T)$, then, via \rf{MIR}, we get
$$\lim_{R}\frac 1R\lim_m\intl{I^m_R}\dm \n v^m(\tau)\dm_2^4d\tau\leq 2\lim_R\intl{I_R}\dm \n v(t)\dm_2^2d\tau=0\,.$$ Letting $m\to\infty$ and then letting $R\to\infty$ in formula \rf{MEE}, we arrive at \rf{EE-IV} with $E(t,s)$ given by \rf{SEE-IV}. Finally, we justify \rf{SUM}. First of all we stress that by virtue of the energy inequality $\dm v^m(t)\dm_2<\dm v^m(s)\dm_2$ true for all $t>s$ and $m\in\N$, we get that in \rf{DIFF} only the last sum can be positive. Concerning \rf{SUM} it is enough to restrict ourselves to the finite sum with indexes $h\in \{1,\dots,p_{12}(m,R)\}\cup\{1,\dots,p_{21}(m,R)\}$. Taking into account that $p_{12}(m,R)=p_{21}(m,R)$, for all $h\in 
\{1,\dots,p_{12}(m,R)\}$  we have $$\ba{l}\dm v^m(t_h^{12})\dm_2^2-\dm v^m(s_h^{12})\dm_2^2<0\,,\VS \dm v^m(t_h^{21})\dm_2^2-\dm v^m(s^{21}_h)\dm_2^2<0\,.\ea$$ Finally, if  $h\in\{2,\dots, p_{12}(m,R)\}$, then the following holds: $$\dm v^m(t_h^{12})\dm_2^2-\dm v^m(s_{h-1}^{21})\dm_2^2<0\,.$$ Hence we obtain 
$$\ba{l}\!\!\mbox{${\underset{h=1}{\overset {p_{21}(m,R)}\sum}}$}\Big[\dm  v^m(t_{h}^{21})\dm_2^2-2\dm v^m(s_h^{21})\dm_2^2\Big] +\!\mbox{${\underset{h=1}{\overset {p_{12}(m,R)}\sum}}$}\Big[2\dm  v^m(t_h^{12})\dm_2^2-\dm v^m(s_h^{12})\dm_2^2\Big]\VS\hskip7cm =B(m,R)+\dm v^m(t_1^{12})\dm_2^2-\dm v^m(s_{p_{21}}^{21})\dm_2^2\,,\ea$$ with $B(m,R)< 0$ and $t_1^{12}<s^{21}_{p_{21}}$\,. 
{\small

\end{document}